\documentclass[11pt,leqno,english,babel]{amsart}
\usepackage{amssymb, amsmath, amsfonts, epsf, epic, eepic}
\usepackage{t1enc, amsthm, amsopn}
\usepackage[latin1]{inputenc}
\usepackage[all]{xy}

\usepackage{color}

\usepackage{amscd}

\newtheorem{Thm}{Theorem}[section]
\newtheorem{Lem}[Thm]{Lemma}
\newtheorem{Def}[Thm]{Definition}
\newtheorem{Cor}[Thm]{Corollary}
\newtheorem{Prop}[Thm]{Proposition}
\newtheorem{Ex}[Thm]{Example}
\newtheorem{Rem}[Thm]{Remark}
\newtheorem{Ques}[Thm]{Question}

\title[Simple-minded systems]{Simple-minded systems in stable module
categories }

\author{Steffen Koenig and Yuming Liu$^*$}

\address{Steffen Koenig
\newline Universit\"at zu K\"oln
\newline Mathematisches Institut
\newline Weyertal 86-90
\newline D-50931 K\"oln
\newline Germany}
\email{skoenig@mi.uni-koeln.de}

\address{Yuming Liu
\newline School of Mathematical Sciences
\newline Laboratory of Mathematics and Complex Systems
\newline Beijing Normal
University
\newline Beijing 100875
\newline P.R.China}
\email{ymliu@bnu.edu.cn}

\date{version of \today}

\newenvironment{Proof}[1][Proof]{\begin{trivlist}
\item[\hskip \labelsep {\bfseries #1}]}{\flushright
$\Box$\end{trivlist}}

\newcommand{\lra}{\longrightarrow}

\newcommand{\ra}{\rightarrow}
\newcommand{\sdp}{\times\kern-.2em\vrule height1.1ex depth-.05ex}
\newcommand{\epi}{\lra \kern-.8em\ra}

\newcommand{\cal}{\mathcal}

 \normalbaselines
\thanks{The corresponding author has been supported by a
Marie Curie Fellowship IIF}

\begin{document}

\renewcommand{\thefootnote}{\alph{footnote}}
\setcounter{footnote}{-1} \footnote{ $^*$ Corresponding author.}
\renewcommand{\thefootnote}{\alph{footnote}}
\setcounter{footnote}{-1}

\renewcommand{\thefootnote}{\alph{footnote}}
\setcounter{footnote}{-1} 

\begin{abstract}
Simple-minded systems in stable module
categories are defined by orthogonality and
generating properties so that the images of the simple modules under a stable
equivalence form such a system. Simple-minded systems are shown to be
invariant under stable equivalences; thus the set of all simple-minded
systems is an invariant of a stable module category. The simple-minded systems
of several classes of algebras are described and connections to the
Auslander-Reiten conjecture are pointed out.

\end{abstract}

\maketitle

\section{Introduction}

Three categories are usually associated with a finite dimensional
algebra $A$: The module category $mod A$, which is an abelian category,
the derived category $D^b(mod A)$, which is triangulated, and the stable
category $\underline{mod}A$, which is also triangulated in case $A$ is
self-injective.

The abelian category $mod A$ is generated by the set of simple
modules and, in a different sense, by each progenerator, that is,
by a full set of indecomposable
projective modules. Equivalences of module categories are described by
Morita theory, in terms of images of projective modules or by progenerators.
The triangulated category $D^b(mod A)$ is also generated by the set of
simple modules and alternatively by each tilting complex. Equivalences of
derived categories are described by Rickard's and Keller's versions of
Morita theory, again in terms of images of projective modules.
Rickard \cite{Rickard2002} has shown how to assign a tilting complex to
a set of objects 'behaving like simple modules', thus allowing to switch
between the two kinds of generators. Rouquier \cite{Rouquier2008}
formalised a concept of
generators of triangulated categories and used it to define the
dimension of a triangulated category.

The stable category $\underline{mod}A$ is generated by the set of
simple modules, too. But the projective modules are not visible in
this category and there is no substitute known for progenerators. An
analogue of Morita theory for stable categories is missing. In
particular, it is not known how to characterize equivalences of
stable categories in terms of images of generators. In fact, it is
not even known how to best define 'generators' of stable categories.
This appears to be a major obstruction to solve a fundamental
problem on stable categories, the Auslander-Reiten conjecture; this
conjecture states that stable equivalences preserve the number of
isomorphisms classes of non-projective simple modules.

\emph{The aim of this article is to suggest and
to explore a new definition of generating sets of stable categories,}
which includes the set of (non-projective) simple modules as an example.

This suggestion is not the first one made. Pogorzaly
\cite{Pogorzaly1993,Pogorzaly1994} introduced what he called maximal systems
of orthogonal stable bricks. He showed that these generate the stable
Grothendieck group, and he used his concept to prove the Auslander-Reiten
conjecture for self-injective special biserial algebras.

The main features of Pogorzaly's systems are mutual orthogonality
and maximality. These properties are clearly invariant under stable
equivalences, while it is a problem to show finiteness of the system
in general. In contrast to this, the generating systems we are
introducing here - the simple-minded systems - satisfy, in addition
to the orthogonality properties, a generating condition that
replaces maximality. Simple-minded systems always are finite
(Proposition \ref{basicfacts}). Invariance under stable equivalences
is not for free any more, but it is true; this is one of the main
results we are going to prove (Theorem \ref{invariant}). The
generating assumption we are using also provides a direct relation
to the stable Grothendieck group.

We consider simple-minded systems for triangular algebras and one-point
extensions and for Nakayama algebras, pointing out connections with the
Auslander-Reiten conjecture. Moreover, we compare the new concept with that
of Pogorzaly and we use the results of this comparison to define the concept
of stable Loewy length. 

\section{Definition and basic properties}

Let $R$ be a commutative artin ring. Recall from \cite{ARS} that an
$R$-algebra $A$ is called an artin algebra if $A$ is finitely
generated as a $R$-module. Important examples of artin algebras are
finite dimensional algebras over a field.

Given an artin algebra $A$, we denote by mod$A$ the category of all
finitely generated left $A$-modules.
For an $A$-module $X$, we denote by $soc(X)$, $top(X)$, and $rad(X)$
its socle, top and radical, respectively. We denote by mod$_{\cal
P}A$ the full subcategory of mod$A$ consisting of modules without
direct summands isomorphic to a projective module. For an $A$-module
$X$, there is a maximal summand (unique up to isomorphism) which has
no nonzero projective summands. We call this summand the
non-projective part of the module $X$.

The stable category $\underline{mod}A$ of $A$ is defined as follows:
The objects of $\underline{mod}A$ are the same as those of mod$A$,
and the morphisms between two objects $X$ and $Y$ are given by the
quotient $R$-module $\underline{\mbox{Hom}}_A(X,Y)$ =
Hom$_A(X,Y)$/$\mathcal{P}(X,Y)$, where $\mathcal{P}(X,Y)$ is the
$R$-submodule of Hom$_A(X,Y)$ consisting of those homomorphisms from
$X$ to $Y$ which factor through a projective $A$-module.

Given two artin algebras $A$ and $B$, we say that $A$ and $B$ are
stably equivalent if their stable categories $\underline{mod}A$ and
$\underline{mod}B$ are equivalent. The Auslander-Reiten translate
$\tau=DTr$ over an artin algebra and the Heller functor (i.e. the
syzygy functor) $\Omega$ over a self-injective algebra are typical
examples of stable self-equivalences. For basic material on
stable equivalence, we refer the reader to \cite{AR1973},
\cite{AR1978}, \cite{ARS}.

Let $A$ be an artin algebra. In \cite{AR1978}, Auslander and Reiten
defined $e(A)$ to be the full additive subcategory of mod$A$ whose
indecomposable objects are the indecomposable non-injective objects
$X$ in mod$A$, such that if $ 0 \longrightarrow X \longrightarrow Y
\longrightarrow Z \longrightarrow 0 $ is an almost split sequence,
then $X$ or $Y$ is projective. Using the notion of node introduced
by Martinez-Villa (cf. \cite{MV1990}), the
indecomposable objects of $e(A)$ consist of precisely the following
three classes of modules: simple projective modules, nodes, and
indecomposable non-simple non-injective projective modules. Clearly
$e(A)$ has only a finite number of indecomposable modules. We denote
by $e'(A)$  the full subcategory of $e(A)$ whose indecomposable
objects are simple projective modules and nodes.

Let $\cal C$ be a class of $A$-modules. We denote by $\langle\cal
C\rangle$ the full subcategory of mod$A$ consisting of modules which
are direct summands of finite direct sums of objects in $\cal C$.
For two subcategories $\cal C$ and $\cal D$ of mod$A$, we denote by
$\langle\cal C\rangle
* \langle\cal D\rangle$ the class of indecomposable $A$-modules $Y$ such that
there is
a short exact sequence of the following form
$$ (\dag) \ \ \ \ \  0
\longrightarrow X \longrightarrow Y\oplus P \longrightarrow Z
\longrightarrow 0, $$ where $Z\in \langle\cal D\rangle, X\in
\langle\cal C\cup e(A)\rangle$, and $P$ is a projective $A$-module.
We put $\langle\cal C\rangle_1=\langle\cal C\rangle$ and we define
inductively $\langle\cal C\rangle_n=\langle\langle\cal
C\rangle_{n-1}*\langle\cal C\rangle\rangle$ for $n\geq 2$.

\begin{Def}
\label{simple-mindedsystem} Let $A$ be an artin algebra. A class
of objects $\mathcal{S}$ in mod$_{\cal P}A$ is called a
\emph{simple-minded system} (for short: s.m.s.) if the following two
conditions are satisfied:

{\rm (1)} {\rm(}orthogonality condition{\rm)} For any $S,T\in
\mathcal{S}$, $$\underline{Hom}_A(S,T)=\left\{\begin{array}{cc} 0 &
S\neq T, \\ \mbox{division ring,} & \mbox{S=T.}
\end{array}\right.$$

{\rm (2)} {\rm(}generating condition{\rm)} For any indecomposable
non-projective $A$-module $X$, there exists some natural number $n$
(depending on $X$) such that $X\in \langle\cal S\rangle_n.$
\end{Def}

\begin{Rem} $(1)$ The definition of a simple-minded system formally
depends on the chosen algebra $A$. In Theorem \ref{invariant} we will see
that in fact a simple-minded system depends only on the equivalence class
of the stable module category $\underline{mod}A$. Therefore, instead of
talking of a simple-minded system over $A$ we may then also talk of a
simple-minded system in $\underline{mod}A$.

$(2)$ By definition, there is no simple-minded system over a
semisimple algebra. From now on, we assume that all algebras
considered are non-semisimple.

$(3)$ Let $A=B\times C$ be a direct product of artin algebras. Then
it is easy to see that the simple-minded systems over $A$ are
exactly of the forms $\mathcal{S}_1\cup \mathcal{S}_2$, where
$\mathcal{S}_1$ is a simple-minded system over $B$ and
$\mathcal{S}_2$ is a simple-minded system over $C$.

$(4)$ When the algebra $A$ is self-injective, its stable module category is
triangulated. In this setup, parallel and independent work of Rickard and 
Rouquier \cite{RR2010} is discussing the problem of reconstructing $A$ from 
its stable module category. They are also using (\cite{RR2010}, 3.2, 
hypothesis 1) the orthogonality and generating conditions satisfied
by the simple modules. Moreover, they are adding a splitting field assumption
and a condition using that there are no extensions between simple modules in
negative degrees. To formulate the latter condition needs the triangulated 
structure. 
\end{Rem}

The following lemma is an easy consequence of our definition.

\begin{Lem} \label{consequence} Let $A$ be an artin algebra.

$(1)$  Suppose that $\mathcal{S}$ is a simple-minded system over
$A$. Then for any $X\in \mathcal{S}$, $X$ is an indecomposable
non-projective module. Moreover, the objects in $\mathcal{S}$ are
pairwise non-isomorphic.

$(2)$ Let $\mathcal{S}$ be a complete set of non-isomorphic simple
non-projective $A$-modules. Then $\mathcal{S}$ is a simple-minded
system over $A$.

$(3)$  If $\mathcal{S}$ is a simple-minded system, then
$\mathcal{S}$ generates the stable Grothendieck group $G_0^{st}(A)$
of $A$.
\end{Lem}

\begin{Proof} (1) is a direct consequence of the orthogonality condition.
For (2), the orthogonality condition is clear. To prove the generating
condition, we use the natural exact sequence $ 0 \longrightarrow
rad(X) \longrightarrow X \longrightarrow top(X) \longrightarrow 0 $
and induction on the Loewy length $ll(X)$ for an indecomposable
module $X$ in mod$_{\cal P}A$. (3) is an easy consequence of the
definition of $G_0^{st}(A)$ (refer to Remark
\ref{stableGrothendieck}) and the generating condition on
$\mathcal{S}$.

\end{Proof}

\begin{Rem} \label{stableGrothendieck}  $(1)$ Let $A$ be an artin algebra.
Recall from \cite{MV1991} that the stable Grothendieck group
$G_0^{\mathrm{st}}(A)$ is by definition the cokernel of the Cartan
map. In other words, there is the following short exact sequence
$$K_0(A)\stackrel{C_A}{\rightarrow} G_0(A)\rightarrow G_0^{\mathrm{st}}(A)\rightarrow 0,$$
where $C_A$ is the Cartan matrix of $A$ and where $K_0(A)$
(respectively, $G_0(A)$) is a free abelian group of finite rank
generated by isomorphism classes of indecomposable projective
modules (respectively, isomorphism classes of simple modules). For
our purpose, we shall use the following equivalent definition of the
stable Grothendieck group (cf. \cite{MV1990}):
$G_0^{\mathrm{st}}(A)$ is the quotient group
$\mathcal{L}_1/\mathcal{R}_1$, where $\mathcal{L}_1$ is the free
group generated by the isomorphism classes $[X]$ of modules $X$ in
mod$_{\cal P}A$, $\mathcal{R}_1$ is the subgroup generated by the
following three classes of elements: $(i)$ $[Y]-[X]-[Z]$, where
$0\longrightarrow X\oplus Q \longrightarrow Y\oplus P
\longrightarrow Z \longrightarrow 0$ is an exact sequence such that
$Q,P$ are projective $A$-modules; $(ii)$ $[X]+[Z]$, where
$0\longrightarrow X\oplus Q \longrightarrow P \longrightarrow Z
\longrightarrow 0$ is an exact sequence such that $Q,P$ are
projective $A$-modules; $(iii)$ $[Z]$, where $0\longrightarrow Q
\longrightarrow P \longrightarrow Z \longrightarrow 0$ is an exact
sequence such that $Q,P$ are projective $A$-modules. In particular,
if $A$ is a self-injective algebra, then $G_0^{\mathrm{st}}(A)$ is
the quotient group $\mathcal{L}_1/\mathcal{R}_1$, where
$\mathcal{L}_1$ is the free group generated by the isomorphism
classes $[X]$ of modules $X$ in mod$_{\cal P}A$, $\mathcal{R}_1$ is
the subgroup generated by the elements $[Y]-[X]-[Z]$ such that there
is an exact sequence $0\longrightarrow X \longrightarrow Y\oplus P
\longrightarrow Z \longrightarrow 0$ with $P$ projective.

$(2)$ The statement (3) in Lemma \ref{consequence} says that the
simple-minded systems are just sets of generators of stable Grothendieck
group. For an artin algebra $A$ of finite global dimension it can be
easily proved that the stable Grothendieck group is trivial, while
on the other hand, there still may exist nontrivial stable
equivalence related to $A$ and such information will be recorded in
the simple-minded systems over $A$ (see Theorem
\ref{invariant}). This indicates that the simple-minded system is a
finer notion than the stable Grothendieck group.
\end{Rem}

Before giving further properties of a simple-minded system, we prove
the following lemma, to be used frequently.

\begin{Lem} \label{leftorthogonal} Let $ 0 \longrightarrow X
\longrightarrow Y \longrightarrow Z \longrightarrow 0 $ be an exact
sequence of $A$-modules and $M$ an $A$-module. If
$\underline{Hom}_A(M,X)=\underline{Hom}_A(M,Z)=0$, then
$\underline{Hom}_A(M,Y)=0$.
\end{Lem}

\begin{Proof} Applying the functor Hom$_A(-,M)$ to the above exact sequence and
using Auslander-Reiten formula, we get the following exact
commutative diagram
$$\begin{CD}   Ext_A^1(Z,M) @>>> Ext_A^1(Y,M)@>>> Ext_A^1(X,M) \\
\wr@VVV \wr@VVV \wr@VVV \\
D\underline{Hom}_A(\tau^{-1}(M),Z) @>>>
D\underline{Hom}_A(\tau^{-1}(M),Y)@>>>
D\underline{Hom}_A(\tau^{-1}(M),X),
\end{CD}$$
where $D=$Hom$_k(-,k)$ denotes the usual duality and $\tau^{-1}=TrD$
is the inverse of the Auslander-Reiten translation. When $M$ runs
through mod$A$, $\tau^{-1}(M)$ runs through mod$_{\cal P}A$. The
lemma thus follows.
\end{Proof}

\begin{Rem} \label{non-symmetry} For self-injective algebras, the same
result holds true for $\underline{Hom}_A(-,M)$ and both results are
special cases of \cite[Lemma 1.4]{EK2000}. However,
Lemma~\ref{leftorthogonal} does not hold for
$\underline{Hom}_A(-,M)$ in general. For example, let $A$ be a path
algebra over a field given by the quiver \unitlength=1.00mm
\special{em:linewidth 0.4pt} \linethickness{0.4pt}
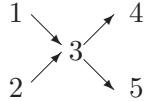
\begin{center}
\begin{picture}(12.00,18.00)
\put(1,12){$1$} \put(9,7){$3$} \put(1,2){$2$} \put(17,12){$4$}
\put(17,2){$5$} \put(4,13){\vector(1,-1){4.0}}
\put(4,4){\vector(1,1){4.0}} \put(11,7){\vector(1,-1){4.0}}
\put(11,9){\vector(1,1){4.0}}
\end{picture}
\end{center}
There is an exact sequence of $A$-modules $ 0 \longrightarrow 4
\longrightarrow \begin{array}{cc} 3 \\ 4
\end{array} \longrightarrow 3 \longrightarrow 0 $. Here,
$\underline{Hom}_A(4,\begin{array}{cc} 3 \\ 4
\end{array})=\underline{Hom}_A(3,\begin{array}{cc} 3 \\ 4
\end{array})=0$, but
$\underline{Hom}_A(\begin{array}{cc} 3 \\ 4
\end{array},\begin{array}{cc} 3 \\ 4
\end{array})\neq 0$.
\end{Rem}

The next proposition collects some elementary facts on a
simple-minded system.

\begin{Prop}\label{basicfacts}
Let $A$ be an artin algebra and let $\mathcal{S}$ be a
simple-minded system over $A$. Then we have the following.

$(1)$  $\mathcal{S}$ contains (up to isomorphism) any simple
non-projective injective module.

$(2)$ $\mathcal{S}$ contains (up to isomorphism) any node.

$(3)$ Assume that $\mathcal{S}_1$ and $\mathcal{S}_2$ are two
classes of objects in mod$_{\cal P}A$ such that
$\mathcal{S}_1\subsetneqq \mathcal{S}\subsetneqq \mathcal{S}_2$.
Then neither $\mathcal{S}_1$ nor $\mathcal{S}_2$ is a
simple-minded system.

$(4)$ The number of objects in $\mathcal{S}$ is finite, that is, the
cardinality $|\mathcal{S}|<\infty$.
\end{Prop}

\begin{Proof} (1) Let $S$ be a simple non-projective injective
module. Suppose that $\mathcal{S}$ does not contain $S$. Then $S$
can be generated by an exact sequence of the following form
$$0 \longrightarrow X
\longrightarrow S\oplus P\longrightarrow Z \longrightarrow 0, $$
where $X\in$ mod$A$, $Z\in \langle\mathcal{S}\rangle$ and $P$ is a
projective $A$-module. Moreover, we can assume that the morphism
$S\longrightarrow Z$ in the above exact sequence is nonzero. But
then $S\longrightarrow Z$ splits and therefore $S$ is a direct
summand of $Z$. This contradicts the assumption that
$\mathcal{S}$ does not contain $S$ and (1) follows.

\medskip
(2) Similarly, suppose that $S$ is a node and that
$\mathcal{S}$ does not contain $S$. Then $S$ can be generated by
an exact sequence of the following form
$$0 \longrightarrow X
\longrightarrow S\oplus P\longrightarrow Z \longrightarrow 0, $$
where $X\in$ mod$A$, $Z\in \langle\mathcal{S}\rangle$ and $P$ is a
projective $A$-module. Moreover, we may assume that neither $X$ nor
$Z$ contains a summand isomorphic to $S$. Since $S$ is a node, we
have an almost split sequence $0 \longrightarrow S \longrightarrow
E\longrightarrow \tau^{-1}(S) \longrightarrow 0$ with $E$
projective. Since $X$ contains no summand isomorphic to $S$, any
homomorphism $S\longrightarrow X$ must factor through the left
almost split homomorphism $S\longrightarrow E$ and therefore
$\underline{Hom}_A(S,X)=0$. Similarly, $\underline{Hom}_A(S,Z)=0$.
It follows from Lemma \ref{leftorthogonal} that
$\underline{Hom}_A(S,S)=0$. This contradiction shows that $S$ is an
object in $\mathcal{S}$.

\medskip
(3) We only need to prove the statement for $\mathcal{S}_1$. Suppose
that $\mathcal{S}_1$ is a simple-minded system and that $X\in
\mathcal{S}\setminus \mathcal{S}_1$. Then $X$ is generated from
objects in $\langle\mathcal{S}_1\cup e(A)\rangle$. By (2) and the
definition of simple-minded system, $X$ is left
orthogonal (in the stable category) to every object in
$\langle\mathcal{S}_1\cup e(A)\rangle$. It follows from Lemma
\ref{leftorthogonal} that $\underline{Hom}_A(X,X)=0$. This is a
contradiction and therefore $\mathcal{S}_1$ is not a simple-minded
system.

\medskip
(4) To generate all simple $A$-modules, we only need a finite
number of objects in $\mathcal{S}$, say, $X_1, \cdots, X_n$. We
can assume that $\{X_1, \cdots, X_n\}$ contains all nodes
(otherwise, we just add the nodes into this set). We shall prove
that $\mathcal{S} = \{X_1, \cdots, X_n\}$. Suppose that there
exists $X\in \mathcal{S}\setminus \{X_1, \cdots, X_n\}$. Since
$\underline{Hom}_A(X,X_i)=0$ for all $1\leq i\leq n$, by Lemma
\ref{leftorthogonal}, we have that $\underline{Hom}_A(X,S)=0$ for
any simple module $S$. This is clearly a contradiction and
therefore our conclusion follows.
\end{Proof}

For self-injective algebras, any projective summand of $X_1$
in the generating sequence $(\dag)$ can be cancelled. Thus we get
the following result.

\medskip
\begin{Cor} Let $A$ be a self-injective algebra and let
$\mathcal{S}$ be a simple-minded system over $A$. For any
indecomposable non-projective $A$-module $X$, there is a projective
$A$-module $P$ and a filtration $$X\oplus P=X_0\supseteq
X_1\supseteq \cdots \supseteq X_m=0$$ with the subquotients in
$\langle\mathcal{S}\rangle$.
\end{Cor}

The above corollary suggests the following definition.

\begin{Def}
\label{filtration} Let $A$ be a self-injective algebra and let
$\mathcal{S}$ be a simple-minded system over $A$. For any
indecomposable non-projective $A$-module $X$, we define $\mu(X)$ as
the minimum integer $m$ such that there is a projective $A$-module
$P$ and a filtration
$$X\oplus P=X_0\supseteq X_1\supseteq \cdots \supseteq X_m=0$$ with
the subquotients in $\langle\mathcal{S}\rangle$. For a general
$A$-module $X$, we define $\mu(X)=\max\mu(Y)$ where $Y$ runs through
all the indecomposable non-projective summands of $X$.
\end{Def}


\begin{Prop}\label{filtration-finite} Let $A$ be a self-injective
algebra and let $\mathcal{S}$ be a simple-minded system over $A$.
Then the values $\mu$ can take are bounded above as follows:
For any $A$-module $X$, there is an inequality
$\mu(X)\leq n_0\cdot ll(A)$, where
$ll(A)$ denotes the usual Loewy length of the regular module $A$ and
$n_0=\mu(A/radA)$.
\end{Prop}

\begin{Proof}
The proof proceeds by induction on the usual Loewy-length $ll(X)$
of $X$. Without loss of generality, we can assume that $X$ is an
indecomposable non-projective $A$-module. If $ll(X)=1$, then $X$ is
a simple module and clearly $\mu(X)\leq n_0$. Now assume that
$ll(X)=n>1$. There is an exact sequence $$0 \longrightarrow Y
\longrightarrow X\longrightarrow Z\longrightarrow 0,$$ where $Y,
Z\in$ mod$_{\cal P}A$, $ll(Y)=n-1$ and $Z$ is semisimple. By
induction, there is a projective module $P_1$ such that $Y\oplus
P_1$ has a $\langle\mathcal{S}\rangle$-filtration of length $\leq
(n-1)n_0$, and there is a projective module $P_2$ such that $Z\oplus
P_2$ has a $\langle\mathcal{S}\rangle$-filtration of length $\leq
n_0$. It follows that the module $X\oplus P_1 \oplus P_2$ has a
$\langle\mathcal{S}\rangle$-filtration of length $\leq (n-1)n_0
+n_0=nn_0\leq n_0\cdot ll(A)$.

\end{Proof}

\section{Invariance under stable equivalences}

In this section, we shall prove that the simple-minded systems are
preserved by any stable equivalence. First of all, we recall some
basic facts on functor categories due to Auslander and Reiten
(cf. \cite{AR1973}, \cite{AR1978}).

Let $A$ be an artin algebra. We denote by mod(mod$A$) the category
of finitely presented contravariant functors $F$ from mod$A$ to
abelian groups. By definition, $F\in$ mod(mod$A$) if and only if
there is a morphism $f: X\longrightarrow Y$ in mod$A$ such that
$F$ is the cokernel of the morphism $$(-, f): (-,
X)\longrightarrow (-, Y),$$ where Hom$_A(-,X)=(-,X)$ and
Hom$_A(-,f)=(-,f)$. Moreover, we denote by $\underline{mod}(modA)$
the full subcategory of mod(mod$A$) whose objects are the functors
which vanishes on projective modules. mod(mod$A$) and
$\underline{mod}(modA)$ have enough projective objects and enough
injective objects. There is a natural functor
$\underline{mod}A\longrightarrow \underline{mod}(modA)$ given by
sending $X$ to $\underline{(-,X)}$, where
$\underline{(-,X)}(Y)=\underline{Hom}_A(Y,X)$, which induces an
equivalence between $\underline{mod}A$ and the full subcategory of
projective objects in $\underline{mod}(modA)$. In particular, we
have that two artin algebras $A$ and $B$ are stably equivalent if
and only if the categories $\underline{mod}(modA)$ and
$\underline{mod}(modA)$ are equivalent. Notice also that the
injective objects in $\underline{mod}(modA)$ are of the form
Ext$_A^1(-,X)$ with $X\in$ mod$A$.

The following lemma extends the result in \cite[Lemma 3.4]{AR1978}.

\begin{Lem} \label{ext} Let $\alpha:
\underline{mod}A\longrightarrow \underline{mod}B$ be a stable
equivalence and $X$ be an indecomposable non-injective $A$-module.
Denote also by $\alpha$ the induced equivalence:
$\underline{mod}(modA)\longrightarrow \underline{mod}(modB)$. Then
there is the following correspondence:

$$ \alpha(Ext_A^1(-,X))\simeq \left\{\begin{array}{cc}
Ext_B^1(-,\alpha(X)) & \mbox{if X is not in e(A),} \\
Ext_B^1(-,Y) \mbox{ for some } Y\in e(B) & \mbox{if X is in e(A).}
\end{array}\right. $$Moreover, if $X$ is in $e'(A)$,
then we also have $Y\in e'(B)$.
\end{Lem}

\begin{Proof} Since $X$ is an indecomposable non-injective $A$-module,
the functor
$Ext_A^1(-,X)$ is an indecomposable injective object in
$\underline{mod}(modA)$. It follows that $\alpha(Ext_A^1(-,X))$ is
an indecomposable injective object in $\underline{mod}(modB)$. If
$X$ is not in $e(A)$, by \cite[Lemma 3.4]{AR1978}, we have that
$\alpha(Ext_A^1(-,X))\simeq Ext_B^1(-,\alpha(X))$ with $\alpha(X)$
not in $e(B)$. If $X\in e(A)$, then $\alpha(Ext_A^1(-,X))\simeq
Ext_B^1(-,Y)$ for some indecomposable non-injective $B$-module $Y$
since every indecomposable injective object in
$\underline{mod}(modB)$ has this form. We claim that $Y\in e(B)$.
Suppose that $Y$ is not in $e(B)$. Again by \cite[Lemma
3.4]{AR1978}, we know that $X=\alpha^{-1}(Y)$ is not in $e(A)$. This
contradiction shows that $Y\in e(B)$.

Now we suppose that $X\in e'(A)$. Then there is an almost split
sequence$$ 0 \longrightarrow X \stackrel{f}{\longrightarrow} P
\stackrel{g}{\longrightarrow} Z \longrightarrow 0$$with $P$ a
projective $A$-module. By \cite[Proposition 2.1]{AR1978}, we have an
exact sequence
$$0\longrightarrow
\underline{(-,Z)}\longrightarrow Ext_A^1(-,X)$$in
$\underline{mod}(modA)$ such that $Ext_A^1(-,X)$ is an injective
envelope of $\underline{(-,Z)}$. We consider the following almost
split sequence$$ 0 \longrightarrow Y \stackrel{f'}{\longrightarrow}
Q \stackrel{g'}{\longrightarrow} \alpha(Z) \longrightarrow 0,
$$where $Y,Q\in$ mod$B$. Observe that $Q$ is a projective $B$-module
and therefore $Y\in e'(B)$. Again by \cite[Proposition 2.1]{AR1978},
there is an exact sequence
$$0\longrightarrow
\underline{(-,\alpha(Z))}\longrightarrow Ext_B^1(-,Y)$$in
$\underline{mod}(modB)$ such that $Ext_B^1(-,Y)$ is an injective
envelope of $\underline{(-,\alpha(Z))}$. Since under the equivalence
$\alpha: \underline{mod}(modA)\longrightarrow
\underline{mod}(modB)$, the functor $\underline{(-,Z)}$ corresponds
to $\underline{(-,\alpha(Z))}$, it follows that
$\alpha(Ext_A^1(-,X))\simeq Ext_B^1(-,Y)$.

\end{Proof}

\begin{Thm} \label{invariant} Let $\alpha:
\underline{mod}A\longrightarrow \underline{mod}B$ be a stable
equivalence and $\mathcal{S}$ be a simple-minded system over $A$.
Then $\alpha(\mathcal{S})$ is a simple-minded system over $B$.
\end{Thm}

\begin{Proof} Obviously, $\alpha(\mathcal{S})$ is a class of objects
in
 mod$_{\cal P}B$ and satisfies the orthogonality condition in
 $\underline{mod}B$. It remains to prove that $\langle\alpha(\mathcal{S})\cup
 e(B)\rangle$ generates any module in mod$_{\cal P}B$.
First we prove the following

\medskip
{\it Claim:} Let $ 0 \longrightarrow X \stackrel{f}{\longrightarrow}
Y\oplus P \stackrel{g}{\longrightarrow} Z \longrightarrow 0 $ be an
exact sequence in mod$A$ which contains no split exact summands,
where $X\in \langle mod_{\cal P}A\cup e(A)\rangle$, $Y$ and $Z$ are
non-zero, $Z\in \langle\cal S\rangle$, $Y\in$ mod$_{\cal P}A$ and
$P$ is a projective $A$-module. Then there is an exact sequence $ 0
\longrightarrow K \stackrel{f'}{\longrightarrow} \alpha(Y)\oplus P'
\stackrel{g'}{\longrightarrow} \alpha(Z) \longrightarrow 0 $ in
mod$B$ which contains no split exact summands, where $K\in
\langle\alpha(X)\cup e(B)\rangle$ and $P'$ is a projective
$B$-module.

\medskip
{\it Proof of Claim:} Since $ 0 \longrightarrow X
\stackrel{f}{\longrightarrow} Y\oplus P
\stackrel{g}{\longrightarrow} Z \longrightarrow 0 $ is an exact
sequence in mod$A$ with no split exact summands, by
\cite[Proposition 2.1]{AR1978}, we know that
$$\underline{(-,Y)}\stackrel{\underline{(-,g)}}{\longrightarrow}
\underline{(-,Z)}\longrightarrow F\longrightarrow 0$$ is a minimal
projective presentation of $F=Coker(-,g)$ in
$\underline{mod}(modA)$, and that
$$0\longrightarrow F\longrightarrow Ext_A^1(-,X)\stackrel{Ext_A^1(-,f)}{\longrightarrow}
Ext_A^1(-,Y)$$ is a minimal injective presentation of
$F=Coker(-,g)$ in $\underline{mod}(modA)$.

Since $\alpha: \underline{mod}A\longrightarrow \underline{mod}B$
is a stable equivalence, we can choose $g'':
\alpha(Y)\longrightarrow \alpha(Z)$ such that
$\alpha(\underline{g})=\underline{g''}$ and choose $t:
P'\longrightarrow \alpha(Z)$ such that $P'\longrightarrow
\alpha(Z)\longrightarrow Cokerg''$ is a projective cover. Let
$g'=(g'',t): \alpha(Y)\oplus P' \longrightarrow \alpha(Z)$, and
consider the exact sequence
$$ 0 \longrightarrow K
\stackrel{f'}{\longrightarrow} \alpha(Y)\oplus P'
\stackrel{g'}{\longrightarrow} \alpha(Z) \longrightarrow 0 .$$
Clearly this sequence has no split exact summands. So, again by
\cite[Proposition 2.1]{AR1978}, there is an exact sequence
$$\underline{(-,\alpha(Y))}\stackrel{\underline{(-,\alpha(g))}}{\longrightarrow}
\underline{(-,\alpha(Z))}\longrightarrow G\longrightarrow
Ext_B^1(-,K)$$ in $\underline{mod}(modB)$, where
$G=Coker(-,\alpha(g))=\alpha(F)$ and $G\longrightarrow
Ext_B^1(-,K)$ is an injective envelope of $G$. It follows that
$\alpha(Ext_A^1(-,X))\simeq Ext_B^1(-,K)$. Write down $X=X_1\oplus
\cdots \oplus X_m \oplus X_{m+1} \oplus\cdots \oplus X_n$, where
each $X_i~(1\leq i\leq n)$
 is an indecomposable non-injective $A$-module
and $X_i$ is not in $e(A)$ for $1\leq i\leq m$, $X_i\in e(A)$ for
$m+1\leq i\leq n$. By Lemma \ref{ext}, $Ext_B^1(-,K)\simeq
\alpha(Ext_A^1(-,X))\simeq
\bigoplus_{i=1}^{n}\alpha(Ext_A^1(-,X_i))=Ext_B^1(-,\alpha(X_1))\oplus
\cdots \oplus Ext_B^1(-,\alpha(X_m))\oplus \cdots \oplus
Ext_B^1(-,K_n)$ for some $K_i~(m+1\leq i\leq n)\in e(B)$. Therefore
$K\simeq \alpha(X_1)\oplus \cdots \oplus \alpha(X_m) \oplus K_{m+1}
\oplus\cdots \oplus K_n\in \langle\alpha(X)\cup e(B)\rangle$. This
finishes the proof of Claim.

From the above claim, it is easy to see that $\alpha(\langle\cal
S\rangle_n)=\langle\alpha(\cal S)\rangle_n$ for any natural number
$n$. It follows that $\langle\alpha(\mathcal{S})\cup
 e(B)\rangle$ generates any module in mod$_{\cal P}B$.

\end{Proof}

The theorem shows that simple-minded systems are stably
invariant. As an application, we determine the simple-minded systems
over the $4$-dimensional weakly symmetric local $k$-algebra
$A_t=k<x,y>/(x^2,y^2,xy-tyx)$, where $k$ is an algebraically closed
field and $t\neq 0$ is an element in $k$. If $chark=2$ and $t=1$,
then $A_t$ is isomorphic to the group algebra of the Klein 4-group.

Let $S$ be the unique (up to isomorphism) simple $A_t$-module.
The Auslander-Reiten quiver of $A_t$ is known to have a component
$\mathcal{C}$ containing $S$ and a $\mathbb{P}_1(k)$-family of
homogenous tubes. For any non-projective module $X$ in the component
$\mathcal{C}$, $X$ is an image of $S$ under some appropriate
composition of the stable equivalence functors $DTr$ and $\Omega$.
It follows that each non-projective $A_t$-module $X\in \mathcal{C}$
defines a simple-minded system over $A_t$. Conversely, every
simple-minded system over $A_t$ is of this form. In fact, if there
is another kind of simple-minded system $\mathcal{S}$ over $A_t$,
then by Proposition \ref{basicfacts} $(3)$, $\mathcal{S}$ should not
contain any module in $\mathcal{C}$. So each $X\in \mathcal{S}$ has
even dimension. Since the unique indecomposable projective
$A_t$-module is also even-dimensional, $\mathcal{S}$ can only
generate even-dimensional modules, a contradiction!

The above argument works for all group algebras $A$ of finite
$p$-groups if the characteristic of the field $k$ is the prime number
$p$. Indeed, by a result of Carlson (\cite{Carlson1998}), an
$A$-module $M$ satisfies $\underline{Hom}_A(M,M)=k$ if and only if
$M$ is an endotrivial module. On the other hand, each endotrivial
module $M$ induces a stable self-equivalence of Morita type over $A$
such that the unique simple module $k$ is mapped to $M$. Therefore every
endotrivial module defines a simple-minded system and these are the
all simple-minded systems over $A$. Theorem \ref{invariant} implies:

\begin{Cor} The Auslander-Reiten conjecture holds true for two algebras
(i.e. two algebras related by a stable equivalence have the same number
of non-projective simple modules up to isomorphism)
if one of the algebras is a group algebra of a finite $p$-group in
characteristic $p$.
\end{Cor}

This result is due to Linckelmann (\cite[Theorem
3.4]{Linckelmann1996}).

\medskip
Clearly, in the above examples, each simple-minded system can be
obtained from the simple modules by applying a suitable stable
self-equivalence. However, the following example shows that
simple-minded systems over an artin algebra are in general not acted
upon transitively by the group of stable self-equivalences. It may
be interesting to determine all the artin algebras with transitive
action of the stable self-equivalences on the simple-minded systems.

\begin{Ex} \label{non-transitivity} Let $k$ be a field. Let $A$ be a finite
dimensional $k$-algebra with the following regular representation
\unitlength=1.00mm \special{em:linewidth 0.4pt}
\linethickness{0.4pt}
\begin{center}
\begin{picture}(12.00,12.00)
\put(-11,3){$A$}\put(-6,3){$=$}\put(0,7){$1'$} \put(0,3){$2'$}
\put(0,-1){$1'$} \put(4,3){$\oplus$} \put(12,7){$2'$}
\put(8,3){$1'$}
\put(16,3){$3'$}\put(12,-1){$2'$}\put(20,3){$\oplus$}\put(24,7){$3'$}
 \put(24,3){$2'$}\put(24,-1){$3'$}
\end{picture}
\end{center}
\vspace{0.1cm} and let $B$ be a finite dimensional $k$-algebra with
the following regular representation \unitlength=1.00mm
\special{em:linewidth 0.4pt} \linethickness{0.4pt}
\begin{center}
\begin{picture}(12.00,14.00)
\put(-11,4){$B$} \put(-6,4){$=$} \put(0,10){$1$} \put(0,6){$2$}
\put(0,2){$3$}\put(0,-2){$1$} \put(4,2){$\oplus$} \put(8,10){$2$}
\put(8,6){$3$}
\put(8,2){$1$}\put(8,-2){$2$}\put(12,2){$\oplus$}\put(16,10){$3$}
\put(16,6){$1$} \put(16,2){$2$}\put(16,-2){$3$}\put(20,-2){$.$}
\end{picture}
\end{center}
\vspace{0.1cm} Both $A$ and $B$ are representation-finite and
symmetric, and there is a stable equivalence of Morita type $\alpha$
between
$B$ and $A$ such that $\alpha(1)=1',\alpha(2)=\begin{array}{c} 3' \\
2'
\end{array},\alpha(3)=\begin{array}{c} 2' \\ 3'
\end{array}$ (cf.
\cite[Section 6]{LiuXi2006}). By Theorem \ref{invariant},
$\{1',\begin{array}{c} 3' \\ 2'
\end{array},\begin{array}{c} 2' \\ 3'
\end{array}\}$ is a simple-minded system. However, there is
no stable self-equivalence $\beta$ over $A$ such that
$\beta(\{1',2',3'\})=\{1',\begin{array}{c} 3' \\ 2'
\end{array},\begin{array}{c} 2' \\ 3'
\end{array}\}$. In fact, if $\beta$ is such a stable
self-equivalence then $\beta$ must be a stable self-equivalence of
Morita type (cf. \cite{Asashiba2003}). Therefore the composition
$\beta^{-1}\alpha$ is a stable equivalence of Morita type between
$B$ and $A$ under which each simple $B$-module corresponds to a
simple $A$-module. It follows from Linckelmann's theorem
(\cite[Theorem 2.1]{Linckelmann1996}) that $B$ and $A$ are Morita
equivalent, which is clearly a contradiction.

\end{Ex}

The next example shows that simple-minded systems over an artin
algebra may even fail to be acted upon transitively by arbitrary stable
equivalences. This is in contrast with the situation for derived
categories. Here, the main result of Rickard's work on derived
equivalences for symmetric algebras in \cite[Theorem
5.1]{Rickard2002} shows transitivity.

\begin{Ex} \label{non-transitivity-2} Let $k$ be a
field. Let $A$ be a finite dimensional $k$-algebra with the
following regular representation \unitlength=1.00mm
\special{em:linewidth 0.4pt} \linethickness{0.4pt}
\begin{center}
\begin{picture}(12.00,18.00)
\put(-11,6){$A$} \put(-6,6){$=$} \put(0,14){$1$} \put(0,10){$2$}
\put(0,6){$1$}\put(0,2){$2$} \put(0,-2){$1$} \put(4,6){$\oplus$}
\put(8,14){$2$} \put(8,10){$1$}
\put(8,6){$2$}\put(8,2){$1$}\put(8,-2){$2$} \put(12,-2){$.$}
\end{picture}
\end{center}
\vspace{0.1cm} Suppose that $B$ is another finite dimensional
$k$-algebra such that there is a stable equivalence $\alpha$ between
$B$ and $A$. As in the above example, $\alpha$ is lifted to a stable
equivalence of Morita type. Since $A$ is an indecomposable
representation-finite symmetric algebra, so is $B$ (by
\cite{Liu2008}). It follows easily that $A$ and $B$ are Morita
equivalent. Without loss of generality we can identify $B$ and $A$
and assume that $\alpha$ is a stable self-equivalence of Morita type
over $A$. By Proposition \ref{Nakayama-algebra}, one can verify that
there are precisely four simple-minded systems over $A$:

$\mathcal{S}_1 = \{1,2\}$; $\mathcal{S}_2 = \{1,\begin{array}{c} 2 \\
 1 \\ 2
\end{array}\}$; $\mathcal{S}_3 = \{\begin{array}{c} 1 \\
 2 \\ 1
\end{array},2\}$; $\mathcal{S}_4 = \{\begin{array}{c} 1 \\
 2 \\ 1 \\ 2
\end{array},\begin{array}{c} 2 \\
 1 \\ 2 \\ 1
\end{array}\}$.\\
 Clearly $\alpha$ commutes with
the syzygy functor $\Omega$, and therefore $\alpha$ cannot map the
simple-minded system $\mathcal{S}_1$ to $\mathcal{S}_2$.

\end{Ex}

\section{Simple-minded systems and triangular algebras}

In this section, we apply simple-minded systems to study triangular
algebras and one-point extension algebras.

For simplicity, throughout this section, we consider (finite
dimensional) quiver algebras of the form $kQ/I$, where $k$ is a
field, $Q$ is a quiver and $I$ is an admissible ideal in $kQ$.
Recall that a quiver algebra $kQ/I$ is said to be a triangular
algebra if there is no oriented cycle in its quiver $Q$.

\begin{Prop} \label{triangular-uniqueness} If $A=kQ/I$ is a
triangular algebra, then $A$ has only one simple-minded system
$\mathcal{S}=\{\mbox{simple non-projective A-modules}\}$.
\end{Prop}

\begin{Proof} Clearly we can assume that the quiver $Q$
contains no isolated vertices. Suppose that $\mathcal{S}$ is a
simple-minded system. Then, by Proposition
\ref{basicfacts}, $\mathcal{S}$ must contains all simple injective
$A$-modules, say, $L_{11},\cdots,L_{1i_1}$ (which correspond to the
source vertices in the quiver $Q$ of $A$). By the orthogonal
condition, every other module in $\mathcal{S}$ have no composition
factor isomorphic to $L_{1j}$ $(1\leq j\leq i_1)$. Let
$\{L_{21},\cdots,L_{2i_2}\}$ be simple $A$-modules which correspond
to such vertices $v_L$ that $v_L$ is not a sink vertex but is next
to a source vertex in the quiver $Q$. Take such a simple $A$-module
$L$ which corresponds to a vertex $v_L$, that is, we are in the
following situation: \unitlength=1.00mm \special{em:linewidth 0.4pt}
\linethickness{0.4pt}
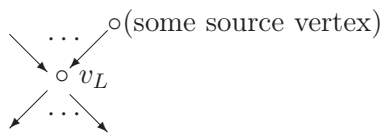
\begin{center}
\begin{picture}(12.00,15.00)
\put(16,12){$\circ$} \put(18,12){$\mbox{(some source vertex)}$}
\put(9,5){$\circ~~v_L$}  \put(8,10){$\cdots$} \put(8,0){$\cdots$}
\put(3,11.5){\vector(1,-1){5.0}} \put(8,4){\vector(-1,-1){5.0}}
\put(11,3.5){\vector(1,-1){5.0}} \put(16,12){\vector(-1,-1){5.0}}
\end{picture}
\end{center}
\vspace{0.1cm} Then $L$ can be generated by an exact sequence of the
following form
$$0 \longrightarrow X
\longrightarrow L\oplus P\longrightarrow Z \longrightarrow 0, $$
where $X\in$ mod$A$, $Z\in \langle\mathcal{S}\rangle$ and $P$ is a
projective $A$-module. Thus there exists some $S\in \mathcal{S}$
such that $L\subseteq soc(S)$. However, if $L_{2j}$ $(1\leq j\leq
i_2)$ is a composition factor of a module in $\mathcal{S}$, then
$L_{2j}$ must occurs in its top. It follows that $L$ is also
contained in $top(S)$. This forces $L\simeq S$ by the
indecomposability property. We thus proved that $\mathcal{S}$
contains the simple modules $L_{21},\cdots,L_{2i_2}$. Observe that
the modules other than the above two classes of simple modules in
$\mathcal{S}$ have no composition factor isomorphic to such simple
modules. Next we consider the simple non-projective $A$-modules
$L_{31},\cdots,L_{3i_3}$ which are ``next'' to simple modules
$L_{21},\cdots,L_{2i_2}$. Continuing by this way, we can prove that
$\mathcal{S}$ contains all the simple non-projective $A$-modules and
therefore $\mathcal{S}=\{\mbox{simple non-projective A-modules}\}$.
\end{Proof}

\begin{Cor} \label{triangular-ARconjecture} If $\alpha:
\underline{mod}B\longrightarrow \underline{mod}A$ is a stable
equivalence such that $A$ is a triangular algebra, then $\alpha$
maps each non-projective simple $B$-module to a non-projective
simple $A$-module, and therefore $A$ and $B$ have the same number of
non-projective simple modules.
\end{Cor}

\begin{Proof} Since $\mathcal{S}=\{\mbox{simple non-projective B-modules}\}$ is a simple-minded system over $B$,
$\alpha(\mathcal{S})$ is a simple-minded system. But
$\mathcal{S'}=\{\mbox{simple non-projective A-modules}\}$ is the
only simple-minded system. It follows that
$\alpha(\mathcal{S})=\mathcal{S'}$, and $A$ and $B$ have the same
number of non-projective simple modules.
\end{Proof}

The next result is a special case of the result in \cite[Theorem
4.3]{DM2007}. Note that here we allow the algebras have nodes. Given
two finite dimensional algebras $A$ and $B$. Recall that $A$ and $B$
are said to be stably equivalent of Morita type if there are two
bimodules $_AM_B$ and $_BN_A$ which are projective as left modules
and as right modules such that we have bimodule isomorphisms:
$$ {}_A M\otimes_B N_A\simeq {}_AA_A\oplus {}_AP_A,\ \ \ {}_B
N\otimes_A M_B\simeq {}_BB_B\oplus {}_BQ_B$$ where ${}_AP_A$ and
${}_BQ_B$ are projective bimodules.


\begin{Prop} \label{triangular-Morita} Let $A$ and $B$ be algebras over
a field $k$. Suppose that $A$ and $B$ have no semisimple summands and
that their maximal semsimple quotient algebras are separable.
If two bimodules
$_{A}M_{B}$ and $_{B}N_{A}$ define a stable equivalence of Morita
type between $A$ and $B$ such that $B$ is a triangular algebra, then
$A$ and $B$ are Morita equivalent.
\end{Prop}

To prove Proposition \ref{triangular-Morita}, we need the following
lemma which is a generalization of Linckelmann's result in
\cite[Theorem 2.1(ii)]{Linckelmann1996} for self-injective algebras.

\begin{Lem} \label{st.M-simple}
Let $A$ and $B$ be two indecomposable nonsimple algebras over a field $k$,
whose maximal semisimple quotient algebras are separable.
If two indecomposable bimodules $_{A}M_{B}$ and
$_{B}N_{A}$ define a stable equivalence of Morita type between $A$
and $B$, then $N\otimes_{A}S$ is an indecomposable $B$-module for
each simple $A$-module.
\end{Lem}

\begin{Proof}  Note that under the assumption of this lemma,
by \cite{DM2007} we can assume that
both ($N\otimes_{A}-$, $M\otimes_{B}-$) and ($M\otimes_{B}-$,
$N\otimes_{A}-$) are adjoint pairs. In particular, $N\otimes_{A}-$
and $M\otimes_{B}-$ maps projective (injective, respectively)
modules to projective (injective, respectively) modules, and $P$ and
$Q$ are projective-injective bimodules. We first state two simple
facts.

{\it Fact 1:} For any indecomposable non-(projective-injective)
$A$-module $X$, $N\otimes_{A}X$ is a non-(projective-injective)
$B$-module; Similarly, we have the result for $M\otimes_{B}-$.

Otherwise, $M\otimes_{B}N\otimes_{A}X\simeq X\oplus P\otimes_{A}X$
is a projective-injective $A$-module and so is $X$. This will be a
contradiction!

{\it Fact 2:} For any indecomposable non-(projective-injective)
$A$-module $X$, suppose that $N\otimes_{A}X\simeq Y\oplus E$ with
$Y$ an indecomposable non-(projective-injective) $B$-module. Then $E$
is a projective-injective $B$-module.

Otherwise, let $E=Z\oplus E'$ with $Z$ an indecomposable
non-(projective-injective) $B$-module. Then
$M\otimes_{B}N\otimes_{A}X\simeq X\oplus P\otimes_{A}X\simeq
M\otimes_{B}Y\oplus M\otimes_{B}Z\oplus M\otimes_{B}E'$. The left
hand side of this equality contains only one indecomposable
non-(projective-injective) summand but the right hand side contains at
least two indecomposable non-(projective-injective) summands. A
contradiction!

Now let $S=Ae/radAe$ be a simple $A$-module. We want to show that
$N\otimes_{A}S$ is an indecomposable $B$-module. There are two cases
to be considered.

{\it Case 1.} $Ae$ is a projective-injective $A$-module. In this
case, $D(Ae)=$ Hom$_k(Ae,k)$ is an indecomposable
projective-injective right $A$-module (or equivalently, an
indecomposable projective-injective $A^{op}$-module). Therefore
$D(Ae)\simeq e'A$ for some primitive idempotent $e'$ in $A$. Note
that $soc(e'A)$ is an ideal in $A$ and that $S\simeq soc(e'A)$ as
left $A$-module. For any indecomposable projective-injective
$B$-module $Bf$, $Bf\otimes_{k}e'A$ is a projective-injective
$B\otimes_{k}A^{op}$-module and $soc(Bf\otimes_{k}e'A)\simeq
soc(Bf)\otimes_{k}soc(e'A)$ by the separability assumption. Since as
a $B\otimes_{k}A^{op}$-module, $N$ has no projective summands,
$soc(Bf\otimes_{k}e'A)N=0$. On the other hand,
$soc(Bf\otimes_{k}e'A)N=soc(Bf)e'A)Nsoc(e'A)\simeq
soc(Bf)(N\otimes_{A}soc(e'A))\simeq soc(Bf)(N\otimes_{A}S)=0$. This
implies that the $B$-module $N\otimes_{A}S$ contains no
projective-injective summands and therefore $N\otimes_{A}S$ is
indecomposable by Fact 1 and 2.

{\it Case 2.} $Ae$ is not an injective $A$-module. Since $P$ is a
projective bimodule, we have a decomposition of the following form:
$P=\bigoplus_{i,j}Ae_i\otimes_{k}e_jA$ where $e_i$'s and $e_j$'s are
some primitive idempotents in $A$. Since $P$ is also an injective
bimodule, each $Ae_i$ and each $e_jA$ are also injective modules. It
follows that $e_jA\otimes_{A}(Ae/radAe)=0$ for each above $e_j$ and
$P\otimes_{A}(Ae/radAe)=0$. This implies that
$M\otimes_{B}N\otimes_{A}(Ae/radAe)\simeq (Ae/radAe)$ and therefore
$N\otimes_{A}(Ae/radAe)$ must be an indecomposable $B$-module.
\end{Proof}

{\it Proof of Proposition \ref{triangular-Morita}.} By
\cite{Liu2008}, we can assume that both $A$ and $B$ are
indecomposable $k$-algebras, and that $_{A}M_{B}$ and $_{B}N_{A}$
are indecomposable bimodules. For any simple $A$-module $S$,
Lemma \ref{st.M-simple} implies that $N\otimes_{A}S$ is an
indecomposable $B$-module. We want to show that $N\otimes_{A}S$ is
simple. There are two cases to be considered.

{\it Case 1.} $S$ is non-projective. In this case, Corollary
\ref{triangular-ARconjecture} implies that $N\otimes_{A}S$ is
simple.

{\it Case 2.} $S$ is simple projective. In this case, by \cite[Lemma
3.1]{Liu2006}, $N\otimes_{A}S$ must contain a simple
projective summand and therefore is also simple.

We have proved that $N\otimes_{A}-$ maps each simple $A$-module to a
simple $B$-module. By the generalization of Linckelmann's theorem
(see \cite[Theorem 1.1]{Liu2003}), the functor $N\otimes_{A}-:$
mod$A\longrightarrow$ mod$B$ gives a Morita equivalence.

\ \ \ \ \ \ \ \ \ \ \ \ \ \ \ \ \ \ \ \ \ \ \ \ \ \ \ \ \ \ \ \ \ \
\ \ \ \ \ \ \ \ \ \ \ \ \ \ \ \ \ \ \ \ \ \ \ \ \ \ \ \ \ \ \ \ \ \
\ \ \ \ \ \ \ \ \ \ \ \ \ \ \ \ \ \ \ \ \ \ \ \ \ \ \ \ \ \ \ \ \ \
\ \ \ \ \ \ \ \ \ \ \ \ \ \ \ \ \ \ \ \ \ \ \ \ $\Box$

Finally, we prove a general fact on simple-minded systems of
one-point extension algebras.

\begin{Prop}\label{one-pointextension} Let $B$ be a finite
dimensional algebra over a field $k$ and let
$A=\left(\begin{array}{cc} B & M \\ 0 & k
\end{array}\right)$ be a one-point extension algebra of $B$ by a $B$-module $M$.

$(1)$ If $\mathcal{S}$ is a simple-minded system over $B$, then
$\mathcal{S'}=\mathcal{S}\cup \{L\}$ is a simple-minded system over
$A$, where $L$ is the simple injective $A$-module with projective
cover $\left(\begin{array}{c} M \\ k
\end{array}\right)$.

$(2)$ Each simple-minded system has the form
$\mathcal{S'}=\mathcal{S}\cup \{L\}$ where $\mathcal{S}$ is a
simple-minded system over $B$ and $L$ is as above.
\end{Prop}

\begin{Proof} (1) There is a canonical algebra
epimorphism $A\longrightarrow B$ given by
$\left(\begin{array}{cc} b & m \\ 0 & x
\end{array}\right)\mapsto b$. So every $B$-module is
automatically an $A$-module by this map.
In particular, $B$ is a projective $A$-module and the embedding
functor $_AB\otimes_B-:$ mod$B\longrightarrow$ mod$A$ induces a
functor $_AB\otimes_B-: \underline{mod}B\longrightarrow
\underline{mod}A$. Note that $_AB\otimes_B-:
\underline{mod}B\longrightarrow \underline{mod}A$ is a fully
faithful functor, and that $e(A)\supseteq e(B)$.

Clearly we have $\underline{Hom}_A(X,Y)=0$ for any $X,Y\in
\mathcal{S}$. Since every $A$-module in $\mathcal{S}$ has no
composition factor isomorphic to $L$, we also have
$\underline{Hom}_A(X,L)=\underline{Hom}_A(L,X)=0$ for any $X\in
\mathcal{S}$. This proves the orthogonality condition for
$\mathcal{S'}$. Now let $Y$ be any indecomposable non-projective
$A$-module. Notice that if $L$ is a composition factor of $Y$, then
$L$ must occurs in the top of $Y$. We consider two cases.

{\it Case 1.} $Y$ has no composition factor isomorphic to $L$. In
this case $Y$ is a $B$-module and can be generated by
$\langle\mathcal{S}\cup e(B)\rangle$. Since $e(B)\subseteq e(A)$, we
know that $Y$ is generated by $\langle\mathcal{S'}\cup e(A)\rangle$.

{\it Case 2.} $Y$ contains composition factor isomorphic to $L$.
We have an exact sequence $$ 0 \longrightarrow X \longrightarrow Y
\longrightarrow L^m \longrightarrow 0, $$where $m$ is a natural
number and $X$ contains no composition factor isomorphic to $L$.
It is readily seen that this case is reduced to Case 1.

\medskip
(2) Suppose that $\mathcal{S'}$ is a simple-minded system.
Then, by Proposition \ref{basicfacts}, $\mathcal{S'}$ must contains
the simple injective module $L$. So $\mathcal{S'}=\mathcal{S}\cup
\{L\}$ with $\mathcal{S}$ a class of objects in mod$_{\cal P}A$. For
any $X\in \mathcal{S}$, $X$ contains no composition factor
isomorphic to $L$: otherwise, $L\in top(X)$ and
$\underline{Hom}_A(X,L)\neq 0$. A contradiction! Therefore
$\mathcal{S}\subseteq$ mod$_{\cal P}B$. We shall prove that
$\mathcal{S}$ is a simple-minded system over $B$. Obviously,
$\mathcal{S}$ satisfies the orthogonality condition in
$\underline{mod}B$ since $\underline{mod}B$ is a full subcategory of
$\underline{mod}A$. Now let $Y$ be an indecomposable $B$-module in
mod$_{\cal P}B$. Then without loss of generality we can assume that
the last exact sequence in mod$A$ which generates $Y$ has the
following form:
$$0 \longrightarrow X \longrightarrow Y\oplus
P\longrightarrow Z\longrightarrow 0, $$ where $m$ is a natural
number, $X\in$ mod$A$, $Z\in \langle\mathcal{S}\rangle$ and $P$ is a
projective $B$-module. It follows that $X$ is a $B$-module and
therefore all the exact sequences involved in generating $Y$ lie in
mod$B$. So $Y$ is generated by $\langle\mathcal{S}\cup e(B)\rangle$,
and $\mathcal{S}$ is a simple-minded system over $B$.

\end{Proof}

\begin{Rem}  $(1)$ Using the above proposition, we
get a simple proof of Proposition
\ref{triangular-uniqueness} as follows: without loss of generality
we assume that $A$ is an indecomposable algebra. Therefore $A$ can
be obtained by a finite number of one-point extensions from a single
point and the conclusion follows immediately from Proposition
\ref{one-pointextension}.

$(2)$ The result in Proposition \ref{one-pointextension} can
not be generalized to triangular matrix algebras, i.e. algebras
of the form $\Lambda=\left(\begin{array}{cc} A & M
\\ 0 & B \end{array}\right)$, where $A$ and $B$ are arbitrary algebras.
For example, let
$k$ be an algebraically closed field. Let $A=B=k[x]/(x^3)$ be two
finite dimensional algebras over $k$ and let $M=k[x]/(x^3)$ be the
natural $A$-$B$-bimodule. Consider the triangular matrix algebra
$\Lambda=\left(\begin{array}{cc} A & M
\\ 0 & B
\end{array}\right)\simeq k(\alpha\circlearrowright 1\stackrel{\beta}{\longleftarrow} 2\circlearrowleft
\gamma)/(\alpha^3,\gamma^3,\alpha\beta-\beta\gamma)$.
 Clearly $(x^2)/(x^3)$ is a simple-minded system which corresponds to the simple $\Lambda$-module $1$,
 and $(x)/(x^3)$ is a simple-minded system over
 $B$ which corresponds to the $\Lambda$-module $\begin{array}{cc} 2 \\
 2
\end{array}$. But $\{1,\begin{array}{cc} 2 \\
 2
\end{array}\}$ is not a simple-minded system over $\Lambda$ since $\begin{array}{cc} 2 \\
 2
\end{array}$ is
not self-orthogonal.
\end{Rem}

\section{Simple-minded systems and self-injective algebras}

In this section, we shall compare the simple-minded systems with
Pogorzaly's maximal systems of stable orthogonal bricks over a
self-injective algebra. We simplify Pogorzaly's definition and drop one
condition used by him to exclude a few trivial cases; in this way we arrive
at 'weakly simple-minded system'. For representation finite self-injective
algebras we show these to coincide with the simple-minded systems defined
before. Thus, for these algebras Pogarzaly's concept essentially coincides
with ours. Once this has been achieved, we introduce the notion of stable
Loewy length for modules in a stable category.

Let $A$ be a self-injective algebra over an algebraically closed
field $k$. Recall from \cite{Pogorzaly1994, Pogorzaly1993} that an
indecomposable $A$-module $X$ in $\underline{mod}A$ is said to be a
stable $A$-brick if its stable endomorphism ring $\underline{End}_A(X)$ is
isomorphic to $k$. A family $\{X_i\}_{i\in I}$ of stable $A$-bricks
is said to be a maximal system of orthogonal stable $A$-bricks if
the following conditions are satisfied:

(1) $\tau(X_i)\ncong X_i$ for any $i\in I$;

(2) $\underline{Hom}_A(X_i,X_j)=0$ for any $i\neq j$;

(3) For any nonzero object $X\in \underline{mod}A$, there exists
some $i\in I$ such that $\underline{Hom}_A(X,X_i)\neq 0$ and there
exists some $j\in I$ such that $\underline{Hom}_A(X_j,X)\neq 0$.

Note that the above definition can be simplified. Indeed, one half
of the assumption in condition (3) is enough: the two conditions
``For any nonzero object $X\in \underline{mod}A$, there exists some
$i\in I$ such that $\underline{Hom}_A(X,X_i)\neq 0$'' and ``For any
nonzero object $X\in \underline{mod}A$, there exists some $j\in I$
such that $\underline{Hom}_A(X_j,X)\neq 0$'' are equivalent. This
can be seen from a general fact on stable categories over a
self-injective algebra proved in \cite{Pogorzaly1993}. The general
fact was presented in the proof of \cite[Proposition
1]{Pogorzaly1993}. For convenience of the reader, we include the
proof here.

\begin{Prop} \label{Pogorzaly} $($\cite[Proof of Proposition 1]{Pogorzaly1993}$)$ Let $A$ be a self-injective artin
algebra. Let $M$ be an indecomposable non-projective $A$-module and
$X$ be any indecomposable $A$-module. If there is a nonzero
homomorphism $\underline{f}: X\longrightarrow M$ in
$\underline{mod}A$, then there is a nonzero homomorphism
$\underline{h}: \tau^{-1}\Omega(M)\longrightarrow X$ such that
$\underline{f}\underline{h}\neq 0$ in $\underline{mod}A$.
\end{Prop}

\begin{Proof} If $X\simeq M$, then there is a nonsplit exact sequence in mod$A$:
$0 \longrightarrow \Omega(M) \longrightarrow P \longrightarrow M
\longrightarrow 0$, where $P \longrightarrow M$ is a projective
cover of $M$. It follows from the Auslander-Reiten formula that
$\underline{Hom}_A(\tau^{-1}\Omega(M),M)\simeq
Ext_A^1(M,\Omega(M))\neq 0$.

Assume now that $X\ncong M$. Consider the following exact
commutative diagram
$$\begin{CD}
0@>>> \Omega(M) @>>> P(M)@>{l}>> M@>>> 0 \\
 @. j@VVV i@VVV 1 @VVV \\
0 @>>> Y @>>> X\oplus P(M)@>{(f,l)}>> M@>>> 0,
\end{CD}$$
where $f: X\longrightarrow M$ is a representative of $\underline{f}$
in mod$A$, $l: P(M)\longrightarrow M$ is a projective cover, $i$ is
a canonical embedding and $j$ is induced from $i$. Applying the snake
lemma we get the following exact sequence of $A$-modules: $$0
\longrightarrow \Omega(M) \stackrel{j}{\longrightarrow} Y
\stackrel{s}{\longrightarrow} X \longrightarrow 0.$$ Note that $j$
is not a split monomorphism since otherwise $s$ is a split
epimorphism and therefore $f(X)=0$, a contradiction! By the property
of almost split sequence, we get the following exact commutative
diagram
$$\begin{CD}
0@>>> \Omega(M) @>>> Z@>{r}>> \tau^{-1}\Omega(M)@>>> 0 \\
 @. 1@VVV t@VVV h @VVV \\
0 @>>> \Omega(M) @>{j}>> Y@>{s}>> M@>>> 0,
\end{CD}$$
where the first row is an almost split sequence, $h$ is induced from
$t$. We first note that $h\neq 0$ since otherwise $r$ will be a
split epimorphism, and this is clearly a contradiction! Next we show
that $\underline{h}\neq 0$. Suppose that this is not the case. We
have the following exact commutative diagram
\begin{center} {\setlength{\unitlength}{0.5cm}
\begin{picture}(30,6)
\put(6,5){$0$}
\put(7,5.2){\vector(1,0){2}}\put(9.5,5){$\Omega(M)$}\put(11.6,5.2){\vector(1,0){2}}\put(14,5){$Z$}\put(12.5,5.4){$r'$}\put(16,5.4){$r$}\put(15,5.2){\vector(1,0){2}}\put(17.5,5){$\tau^{-1}\Omega(M)$}
\put(21,5.2){\vector(1,0){2}}\put(23.5,5){$0$}\put(6,5){$0$}
\put(6,2){$0$}\put(7,2.2){\vector(1,0){2}}\put(9.5,2){$\Omega(M)$}\put(12.5,2.5){$j$}\put(11.6,2.2){\vector(1,0){2}}\put(14,2){$Y$}\put(16,1.7){$s$}\put(15,2.2){\vector(1,0){2}}\put(18.7,2){$X$}
\put(21,2.2){\vector(1,0){2}}\put(23.5,2){$0,$}\put(10,3.5){$1$}\put(10.5,4.5){\vector(0,-1){1.5}}\put(13.8,3.5){$t$}\put(14.2,4.5){\vector(0,-1){1.5}}\put(19,4.5){\vector(0,-1){1.5}}
\put(19.2,3.5){$h$}\put(16,3.5){$P$}\put(17.6,4.8){\vector(-1,-1){1.0}}\put(15.8,3.5){\vector(-1,-1){1.0}}\put(17,3.5){\vector(1,-1){1.0}}\end{picture}
}\end{center} where $h$ factors through the projective cover $P$ of
$X$, the homomorphism $P\longrightarrow Y$ is induced from the
property of projective modules, and we denote the composition map
$\tau^{-1}\Omega(M)\longrightarrow P\longrightarrow Y$ by $u$. We
have that $Im(t-ur)\subseteq j(\Omega(M))$ and that
$(t-ur)r'(\Omega(M))=tr'(\Omega(M))=j(\Omega(M))\simeq \Omega(M)$.
Therefore $r'$ is a split monomorphism which is clearly a
contradiction. Finally, let us prove that
$\underline{f}\underline{h}\neq 0$. Indeed, if $fh$ factors through
$P(M)$ then $fh=lh'$ for some $h': \tau^{-1}\Omega(M)\longrightarrow
P(M)$ and consequently $h$ factors through $Y$ . Hence there is
$h_1: \tau^{-1}\Omega(M)\longrightarrow Y$ with $h=sh_1$. Thus
$st=sh_1r$ and $Im(t-h_1r)\subseteq Im(j)$. Then, as
before, a contradiction can be deduced. This shows that
$\underline{f}\underline{h}\neq 0$.
\end{Proof}

\begin{Cor} \label{half} In the definition of maximal system of orthogonal bricks, the two conditions
``For any nonzero object $X\in \underline{mod}A$, there exists some $i\in I$ such that
$\underline{Hom}_A(X,X_i)\neq 0$'' and ``For any nonzero object
$X\in \underline{mod}A$, there exists some $j\in I$ such that
$\underline{Hom}_A(X_j,X)\neq 0$'' are equivalent.
\end{Cor}

\begin{Proof} It suffices to prove it for $X$ indecomposable. Suppose that the condition
``For any nonzero object $X\in \underline{mod}A$, there exists some $i\in I$ such that
$\underline{Hom}_A(X,X_i)\neq 0$'' is satisfied. By Proposition
\ref{Pogorzaly}, $\underline{Hom}_A(X_i,\Omega^{-1}\tau(X))\simeq
\underline{Hom}_A(\tau^{-1}\Omega(X_i),X)\neq 0$. When $X$ runs
through the nonzero objects in $\underline{mod}A$, so does
$\Omega^{-1}\tau(X)$. Therefore we have proved the another
condition. The proof of the other direction is similar.
\end{Proof}

In order to compare simple-minded systems with Pogorzaly's
maximal systems of orthogonal bricks over a self-injective algebra,
we introduce the following definition (note that our definition here
applies in any artin algebra).

\begin{Def}
\label{weaklysimple-mindedsystem} Let $A$ be an artin algebra. A
class of objects $\mathcal{S}$ in mod$_{\cal P}A$ is called a
\emph{weakly simple-minded system} if the following two conditions are
satisfied:

{\rm (1)} {\rm(}orthogonality condition{\rm)} For any $S,T\in
\mathcal{S}$, $$\underline{Hom}_A(S,T)=\left\{\begin{array}{cc} 0 &
S\neq T, \\ \mbox{division ring,} & \mbox{S=T.}
\end{array}\right.$$

{\rm (2)} {\rm(}weak generating condition{\rm)} For any
indecomposable non-projective $A$-module $X$, there exists some
$S\in \mathcal{S}$ (depends on $X$) such that
$\underline{Hom}_A(X,S)\neq 0.$
\end{Def}

\begin{Rem} According to Remark \ref{non-symmetry}, for general artin algebras, the weak generating condition in Definition \ref{weaklysimple-mindedsystem}
is not symmetric, that is, "$\underline{Hom}_A(X,S)\neq 0$" cannot
be replaced by "$\underline{Hom}_A(S,X)\neq 0$".
\end{Rem}

It is easy to see that every simple-minded system is a weakly
simple-minded system. The reason is as follows: Let $A$ be an artin
algebra and let $\mathcal{S}$ be a simple-minded system. To show
that $\mathcal{S}$ is a weakly simple-minded system, we only need to
prove the weak generating condition. Let $0\neq X\in$ mod$_{\cal
P}A$. Suppose that $\underline{Hom}_A(X,T)=0$ for all $T\in
\mathcal{S}$. Then we have that $\underline{Hom}_A(X,S)=0$ for any
simple module $S$ (cf. the proof of Lemma \ref{leftorthogonal}).
This is clearly a contradiction and therefore $\mathcal{S}$
satisfies the weak generating condition. Thus the question arises:
Is every weakly simple-minded system also a simple-minded system?

At least for representation-finite self-injective finite dimensional
algebras, we can prove that the above question has a positive
answer. First we need a lemma. Let $A$ be a finite dimensional
algebra over a field $k$ and let $\mathcal{S} = \{M_1, \cdots,
M_n\}$ be a weakly simple-minded system. Let $X$ be an
$A$-module in mod$_{\cal P}A$. Suppose that
dim$_k\underline{Hom}_A(X,M_i)=d_i$ for $1\leq i \leq n$. Following
\cite{Pogorzaly1993}, we will say that
$\bigoplus_{i=1}^{n}M_i^{d_i}$ is an s-top of $X$ with respect to
$\mathcal{S}$. Of course, s-top($X$) is well-defined for $X$. We
consider the following exact sequence in mod$A$:
$$(*)\ \ \ \ \  0 \longrightarrow X_1 \stackrel{h=(h',h'')}{\longrightarrow} X\oplus P
\stackrel{(f,g)}{\longrightarrow} \mbox{s-top}(X) \longrightarrow 0,
$$
where $f: X\longrightarrow \mbox{s-top}(X)$ is such a morphism that
the coordinates of $\underline{f}$ form a basis of the nonzero
$k$-space $\underline{Hom}_A(X,\mbox{s-top}(X))$ and $g:
P\longrightarrow \mbox{s-top}(X)$ is such a morphism that
$P\longrightarrow \mbox{s-top}(X)\longrightarrow Coker(f)$ is a
projective cover.

\begin{Lem} \label{well-defined}  Let $X$ be an
$A$-module in mod$_{\cal P}A$. Up to isomorphism, the non-projective
part of the module $X_1$ in the above sequence $(*)$ is independent
of the choice of the homomorphism $f: X\longrightarrow
\mbox{s-top}(X)$.
\end{Lem}

\begin{Proof}
First we note that if we replace $g: P\longrightarrow
\mbox{s-top}(X)$ in the above sequence $(*)$ by the projective cover
$g': Q\longrightarrow \mbox{s-top}(X)$, then ker$(f,g')$ and $X_1$
have the isomorphic non-projective part. Now we choose another
homomorphism $f': X\longrightarrow \mbox{s-top}(X)$ such that the
coordinates of $\underline{f'}$ still form a $k$-basis of
$\underline{Hom}_A(X,\mbox{s-top}(X))$. There clearly
is an $A$-module isomorphism $\alpha: \mbox{s-top}(X)\longrightarrow
\mbox{s-top}(X)$ such that $\alpha f-f'$ factors through the
projective cover $q_2: P'\longrightarrow \mbox{s-top}(X)$. More
precisely, there is a homomorphism $q_1: X\longrightarrow P'$ such
that $\alpha f-f'=q_2q_1$. Hence we get the following exact
commutative diagram
$$\begin{CD}
0@>>> Y_1 @>>> X\oplus P'@>{(f,\alpha^{-1}q_2)}>> \mbox{s-top}(X)@>>> 0 \\
 @. d@VVV {\small \left(\begin{array}{cc} 1 & 0
\\ q_1 & 1
\end{array}\right)}@VVV \alpha @VVV \\
0@>>> Y_1' @>>> X\oplus P'@>{(f',q_2)}>> \mbox{s-top}(X)@>>> 0,
\end{CD}$$
where $d$ is induced from the isomorphism $\left(\begin{array}{cc} 1
& 0
\\ q_1 & 1
\end{array}\right)$. It follows that $d$ is an isomorphism.
In particular, $Y_1$ and $Y_1'$ have the isomorphic
non-projective parts.
\end{Proof}

Clearly, if $A$ is a self-injective algebra, $X_1$ contains no
projective summands. However, in general $X_1$ may contain
projective summands (although by our assumption, $X$ contains no
projective summands). According to \cite{Pogorzaly1993}, we define
the s-radical of $X$ with respect to $\mathcal{S}$ to be the
non-projective part of $X_1$ in the above sequence $(*)$. This is
well-defined up to isomorphism, and we shall denote it by
$\mbox{s-rad}(X)$. Moreover, we denote
$\mbox{s-rad}(\mbox{s-rad}^{i-1}(X))$ by $\mbox{s-rad}^i(X)$.

\begin{Thm} \label{equivalence} Let $A$ be a
representation-finite self-injective finite dimensional algebra over
a field $k$ and let $\mathcal{S} = \{M_1, \cdots, M_n\}$ be a weakly
simple-minded system. Then $\mathcal{S}$ even is a
simple-minded system.
\end{Thm}

\begin{Proof}
We only need to prove the generating condition. Let $X$ be an
indecomposable non-projective $A$-module. Suppose that
dim$_k\underline{Hom}_A(X,M_i)=d_i$ for $1\leq i \leq n$. As before,
we consider the following exact sequence in mod$A$:
$$(*)\ \ \ \ \  0 \longrightarrow X_1 \stackrel{h=(h',h'')}{\longrightarrow} X\oplus P
\stackrel{(f,g)}{\longrightarrow} \mbox{s-top}(X) \longrightarrow 0,
$$
where $f: X\longrightarrow \mbox{s-top}(X)$ is such a morphism that
the coordinates of $\underline{f}$ form a basis of the nonzero
$k$-space $\underline{Hom}_A(X,\mbox{s-top}(X))$ and $g:
P\longrightarrow \mbox{s-top}(X)$ is such a morphism that
$P\longrightarrow \mbox{s-top}(X)\longrightarrow Coker(f)$ is a
projective cover. Let $M=\bigoplus_{i=1}^{n}M_i$. Since
$\underline{mod}A$ is a triangulated category (with translation
functor
$\Omega^{-1}:\underline{mod}A\longrightarrow\underline{mod}A$) and
the above exact sequence induces a triangle
$$ X_1 \stackrel{\underline{h}}{\longrightarrow} X\stackrel{\underline{f}}{\longrightarrow} \mbox{s-top}(X)
\stackrel{\underline{e}}{\longrightarrow}\Omega^{-1}(X_1)$$ in
$\underline{mod}A$, after applying the contravariant cohomological
functor $\underline{Hom}_A(-,M)$ to the above triangle, we get the
following long exact sequence of $k$-spaces
$$\cdots \longrightarrow \underline{_(\Omega^{-1}(X),M)}
\stackrel{\underline{(\Omega^{-1}(h),M)}}{\longrightarrow}
\underline{_(\Omega^{-1}(X_1),M)}
\stackrel{\underline{(e,M)}}{\longrightarrow}
\underline{_(\mbox{s-top}(X),M)}
\stackrel{\underline{(f,M)}}{\longrightarrow}
\underline{_(X,M)}\longrightarrow \cdots.
$$
We claim that $\underline{(f,M)}$ is an isomorphism. Indeed, the
spaces $\underline{_(\mbox{s-top}(X),M)}$ and $\underline{_(X,M)}$
have the same $k$-dimension $\Sigma_{i=1}^{n}d_i$ and the canonical
basis elements of $\underline{_(\mbox{s-top}(X),M)}$ map to the
coordinates of $\underline{f}$ which form a basis of
$\underline{_(X,M)}$. It follows that
$\underline{(\Omega^{-1}(h),M)}$ is an epimorphism and that
dim$_k\underline{Hom}_A(\Omega^{-1}(X),M)\geq$
dim$_k\underline{Hom}_A(\Omega^{-1}(X_1),M)$. We can assume that
$X_1\neq 0$ since otherwise $X\simeq \mbox{s-top}(X)\in \langle\cal
S\rangle$ and we are done. Note also that $X_1$ contains no
projective summand. For any indecomposable summand of $X_1$ (we
still denote it by $X_1$), we can similarly take an exact sequence
as $(*)$ in mod$A$:
$$ 0 \longrightarrow X_2 \stackrel{h_1}{\longrightarrow} X_1\oplus
P_1 \stackrel{(f_1,g_1)}{\longrightarrow} \mbox{s-top}(X_1)
\longrightarrow 0.$$From this we also deduce a canonical epimorphism
$\underline{(\Omega^{-1}(h_1),M)}:\underline{_(\Omega^{-1}(X_1),M)}
\longrightarrow \underline{_(\Omega^{-1}(X_2),M)}$ and get an
inequality dim$_k\underline{Hom}_A(\Omega^{-1}(X_1),M)\geq$
dim$_k\underline{Hom}_A(\Omega^{-1}(X_2),M)$. Continuing in this
way, we obtain a sequence of epimorphisms between $k$-spaces:
$$\underline{_(\Omega^{-1}(X),M)}
\stackrel{\underline{(\Omega^{-1}(h),M)}}{\longrightarrow}
\underline{_(\Omega^{-1}(X_1),M)}\stackrel{\underline{(\Omega^{-1}(h_1),M)}}{\longrightarrow}
\underline{_(\Omega^{-1}(X_2),M)}\stackrel{\underline{(\Omega^{-1}(h_2),M)}}{\longrightarrow}
\underline{_(\Omega^{-1}(X_3),M)}\longrightarrow \cdots.
$$The above sequence is induced from the following
sequence
$$\cdots\longrightarrow X_3
\stackrel{\underline{h_2}}{\longrightarrow} X_2
\stackrel{\underline{h_1}}{\longrightarrow} X_1
\stackrel{\underline{h}}{\longrightarrow} X
$$ in $\underline{mod}A$
and the latter one is again induced from the following sequence
$$(**)\ \ \ \ \  \cdots\longrightarrow X_3
\stackrel{h_2'}{\longrightarrow} X_2
\stackrel{h_1'}{\longrightarrow} X_1 \stackrel{h'}{\longrightarrow}
X
$$ in mod$A$. To finish our proof, it suffices to prove the following conclusion:
there exists some natural number $m$ such that $X_m=0$ (and
consequently $X_i=0$ for all $i\geq m$).

By our assumption, all the modules in the above sequence $(**)$ are
indecomposable. We claim that all homomorphisms in $(**)$ are
non-isomorphisms. In fact, if in the original sequence $(*)$ the
s-radical $X_1$ contains an indecomposable summand $X_1'$ such that
$h': X_1'\longrightarrow X$ is an isomorphism, then the inequality
dim$_k\underline{Hom}_A(\Omega^{-1}(X),M)\geq$
dim$_k\underline{Hom}_A(\Omega^{-1}(X_1),M)$ implies that $X_1$ can
not contains any other summands, and therefore $X_1$ must be
isomorphic to $X$. This would leads to the absurd conclusion that
the sequence $(*)$ splits and that $X\simeq X\oplus
\mbox{s-top}(X)$. Similarly, one can show that all $h_i' (i\geq 1)$
are non-isomorphisms. Since $A$ is representation-finite and the
modules in mod$A$ have bounded length, by \cite[Corollary 1.3]{ARS},
for some large $m$ ($m\leq 2^b$, where $b$ denotes the least upper
bound of the lengths of the indecomposable modules in mod$A$), the
composition $h'h_1'\cdots h_m'$ is zero in mod$A$. It follows that
the composition $\underline{(\Omega^{-1}(h_m),M)}\cdots
\underline{(\Omega^{-1}(h_1),M)}\underline{(\Omega^{-1}(h),M)}$ is
zero. Since all $\underline{(\Omega^{-1}(h_i),M)}$ are epimorphisms,
we know that $\underline{Hom}_A(\Omega^{-1}(X_m),M)=0$. By the
weak generating condition, we know that $\Omega^{-1}(X_m)=0$. It
follows that $X_m=0$ since $\Omega^{-1}:
\underline{mod}A\longrightarrow \underline{mod}A$ is an equivalence.
\end{Proof}

\begin{Rem}
\label{acyclic} Suppose that $A$ is any (not necessarily
representation-finite) self-injective algebra over a field $k$.
The above proof implies that for any indecomposable non-projective
$A$-module $X$, $\mbox{s-rad}(X)$ cannot contain a direct summand
isomorphic to $X$. Indeed, if this is the case, we can take all
$X_i$ equals to $X$ in the above proof, and finally we get that
$X=0$, which is a contradiction. Moreover, it is easy to see that
all the modules in $\{\mbox{s-rad}^i(X)|i=0,1,2,\cdots \}$ are
pairwise disjoint, i.e. do not have isomorphic direct summands.
\end{Rem}

It is well-known that the Loewy length is a very useful concept in
the module category mod$A$. It would be interesting to generalize
this notion to the stable module category $\underline{mod}A$. Here,
we replace the simple modules by a simple-minded system. Indeed,
Lemma \ref{well-defined} supplies a way to define the stable
Loewy length of an object in $\underline{mod}A$.

\begin{Def}
\label{stableLoewylength} Let $A$ be a finite dimensional algebra
over a field $k$ and let $\mathcal{S}$ be a simple-minded system
over $A$. For any indecomposable non-projective $A$-module $X$ in
mod$_{\cal P}A$, we define the stable Loewy length of $X$ with
respect to $\mathcal{S}$ (which we denote by s-ll($X$)) to be the
least number $m$ such that $\mbox{s-rad}^m(X)=0$. If there is no
such $m$, then we define $\mbox{s-ll}(X)=\infty$. For a general
module $X\in$ mod$A$, we define $\mbox{s-ll}(X)$ to be the stable
Loewy length of its non-projective part.
\end{Def}

\begin{Cor} \label{finite} Let $A$ be a representation-finite self-injective
finite dimensional algebra over a field $k$ and let
$\mathcal{S}$ be any simple-minded system. Then the stable
Loewy length satisfies the inequality
$\mbox{s-ll}(X)\leq 2^b$ for any $X\in$ mod$A$, where
$b$ denotes the least upper bound of the lengths of the
indecomposable modules in mod$A$.
\end{Cor}

\begin{Proof} This is an easy consequence of the proof of Theorem
\ref{equivalence}.
\end{Proof}

\begin{Ex} Consider the algebra $A$ in Example
\ref{non-transitivity-2}. Both $\mathcal{S}_1 = \{1,2\}$ and $\mathcal{S}_2 = \{1,\begin{array}{c} 2 \\
1 \\ 2
\end{array}\}$ are simple-minded systems over $A$. For any indecomposable non-projective $A$-module
$X$, the stable Loewy length $\mbox{s-ll}(X)$ with respect to
$\mathcal{S}_1$ is equal to the usual Loewy length $ll(X)$. However,
the stable Loewy length $\mbox{s-ll}(X)$ with respect to
$\mathcal{S}_2$ is usually different from $ll(X)$.
For example, the stable Loewy length of the $A$-module $\begin{array}{c} 1 \\
2 \\ 1 \\
2
\end{array}$ with respect to $\mathcal{S}_2$ is equal to $2$ while its usual
Loewy length is $4$.
\end{Ex}

\section{Nakayama algebras}

One motivation to define simple-minded systems is to explore the potential
use of this concept for the Auslander-Reiten conjecture. This conjecture
says that two stably equivalent artin algebras have the same number
of non-isomorphic non-projective simple modules. Based on the
observation in Theorem \ref{invariant}, we pose the following
question.

\begin{Ques} \label{cardinality}
Is the cardinality of each simple-minded system over an artin
algebra $A$ equal to the number of non-isomorphic non-projective
simple $A$-modules?
\end{Ques}

A positive answer to this question implies the Auslander-Reiten
conjecture. Actually, Pogorzaly \cite{Pogorzaly1994} used an
analogous idea in his setup to prove the conjecture for
self-injective special biserial algebras. We think that besides the
relationship with Auslander-Reiten conjecture, Question
\ref{cardinality} is interesting in itself. Proposition
\ref{triangular-uniqueness} shows that the answer is positive for
triangular algebras. The next proposition answers this question for
Nakayama algebras.

\begin{Prop}\label{Nakayama-algebra}
Let $A$ be a Nakayama algebra and let $\mathcal{S}$ be a
simple-minded system. Then the cardinality of $\mathcal{S}$
is equal to the number of non-isomorphic non-projective simple
$A$-modules. Moreover, if we assume that $\mathcal{S} = \{M_1,
\cdots, M_n\}$ and that $\{S_1, \cdots, S_n\}$ is a complete set of
non-isomorphic non-projective simple $A$-modules, then both the set of
tops $top(M_1), \cdots, top(M_n)$ and the set of socles $soc(M_1),
\cdots, soc(M_n)$ coincide, up to ordering, with the set of
simple modules $S_1, \cdots, S_n$.
\end{Prop}

\begin{Proof} First we remind the reader that every
indecomposable module over a Nakayama algebra
is uniserial. Let $S$ be any non-projective simple $A$-module. Then
there exists some $M_i\in \mathcal{S}$ such that $S\simeq soc(M_i)$
by the weak generating condition. This shows that each
non-projective simple $A$-module occurs as a socle of some $M_i\in
\mathcal{S}$. On the other hand, any two different $M_i$ and $M_j$
must have non-isomorphic socles. Indeed, if $M_i$ and $M_j$ satisfy
$soc(M_i)\simeq soc(M_j)$, then there is a monomorphism from one
module to another module, say, $M_i\hookrightarrow M_j$. But clearly
in this case this homomorphism does not factor through a projective
module and therefore $M_i\simeq M_j$ by the orthogonality condition. We
have proved that the cardinality of $\mathcal{S}$ is equal to the
number of non-isomorphic non-projective simple $A$-modules and that
the socle series $soc(M_1), \cdots, soc(M_n)$ is a rearrangement of
$S_1, \cdots, S_n$. To prove the statement for top series, it
suffices to show that any two different $M_i$ and $M_j$ must have
non-isomorphic tops. In fact, if $M_i$ and $M_j$ satisfy
$top(M_i)\simeq top(M_j)$, then there is an epimorphism from one
module to another module, say, $M_i\twoheadrightarrow M_j$. But
clearly this homomorphism does not factor through a projective
module and therefore $M_i\simeq M_j$ by the orthogonality condition.
\end{Proof}

We now give an example to illustrate the above proposition.

\begin{Ex} We consider the Nakayama algebra $B$ in Example
\ref{non-transitivity}. First we display the Auslander-Reiten quiver
of $B$ as follows:

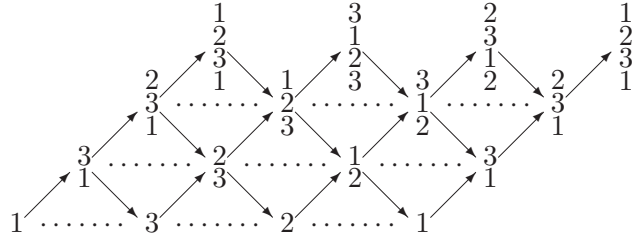
\unitlength=1.00mm \special{em:linewidth 0.4pt}
\linethickness{0.4pt}
\begin{center}
\begin{picture}(40.00,30.00)
\put(-16,0){\vector(1,1){6.0}}\put(-18,-2){$1$} \put(-9,7){$3$}
\put(-9,4){$1$} \put(-7,8){\vector(1,1){6.0}} \put(0,17){$2$}
\put(0,14){$3$}
\put(0,11){$1$}\put(2,16){\vector(1,1){6.0}}\put(9,26){$1$}
\put(9,23){$2$} \put(9,20){$3$}\put(9,17){$1$}
\put(-7,6){\vector(1,-1){6.0}}\put(2,14){\vector(1,-1){6.0}}
\put(11,22){\vector(1,-1){6.0}}
\put(2,0){\vector(1,1){6.0}}\put(0,-2){$3$} \put(9,7){$2$}
\put(9,4){$3$} \put(11,8){\vector(1,1){6.0}} \put(18,17){$1$}
\put(18,14){$2$}
\put(18,11){$3$}\put(20,16){\vector(1,1){6.0}}\put(27,26){$3$}
\put(27,23){$1$} \put(27,20){$2$}\put(27,17){$3$}
\put(11,6){\vector(1,-1){6.0}}\put(20,14){\vector(1,-1){6.0}}
\put(29,22){\vector(1,-1){6.0}}
\put(20,0){\vector(1,1){6.0}}\put(18,-2){$2$} \put(27,7){$1$}
\put(27,4){$2$} \put(29,8){\vector(1,1){6.0}} \put(36,17){$3$}
\put(36,14){$1$}
\put(36,11){$2$}\put(38,16){\vector(1,1){6.0}}\put(45,26){$2$}
\put(45,23){$3$} \put(45,20){$1$}\put(45,17){$2$}
\put(29,6){\vector(1,-1){6.0}}\put(38,14){\vector(1,-1){6.0}}
\put(47,22){\vector(1,-1){6.0}}
\put(38,0){\vector(1,1){6.0}}\put(36,-2){$1$} \put(45,7){$3$}
\put(45,4){$1$} \put(47,8){\vector(1,1){6.0}} \put(54,17){$2$}
\put(54,14){$3$}
\put(54,11){$1$}\put(56,16){\vector(1,1){6.0}}\put(63,26){$1$}
\put(63,23){$2$} \put(63,20){$3$}\put(63,17){$1$}
\put(-14,-2){$\cdots\cdots\cdot$}\put(4,-2){$\cdots\cdots\cdot$}
\put(22,-2){$\cdots\cdots\cdot$}\put(-5,6){$\cdots\cdots\cdot$}
\put(13,6){$\cdots\cdots\cdot$}\put(31,6){$\cdots\cdots\cdot$}
\put(4,14){$\cdots\cdots\cdot$}\put(22,14){$\cdots\cdots\cdot$}
\put(40,14){$\cdots\cdots\cdot$}
\end{picture}
\end{center}
\vspace{0.2cm}where the dotted lines indicate the Auslander-Reiten
translation, and the same vertices are identified. Clearly each
indecomposable non-projective $B$-module is self-orthogonal in
$\underline{mod}B$. Using Proposition \ref{Nakayama-algebra}, it is
not hard to verify that there are precisely five simple-minded
systems over $B$:

$\mathcal{S}_1 = \{1,2,3\}$; $\mathcal{S}_2 = \{1,\begin{array}{cc} 2 \\
 3
\end{array},\begin{array}{ccc} 3 \\
 1 \\ 2
\end{array}\}$; $\mathcal{S}_3 = \{\begin{array}{cc} 1 \\
 2
\end{array},\begin{array}{ccc} 2 \\
 3 \\ 1
\end{array},3\}$; $\mathcal{S}_4 = \{\begin{array}{ccc} 1 \\
 2 \\ 3
\end{array},2,\begin{array}{cc} 3 \\
 1
\end{array}\}$; $\mathcal{S}_5 = \{\begin{array}{ccc} 1 \\
 2 \\ 3
\end{array},\begin{array}{ccc} 2 \\
 3 \\ 1
\end{array},\begin{array}{ccc} 3 \\
 1 \\ 2
\end{array}\}$.\\
On the other hand, since the algebra $A$ in Example
\ref{non-transitivity} is stably equivalent to $B$, there are also
five simple-minded systems over $A$. However, if we consider the
following quotient algebra (which is still a Nakayama algebra but
not self-injective) of $B$ \unitlength=1.00mm \special{em:linewidth
0.4pt} \linethickness{0.4pt}
\begin{center}
\begin{picture}(12.00,14.00)
\put(-11,4){$B'$} \put(-6,4){$=$} \put(0,6){$1$}
\put(0,2){$2$}\put(0,-2){$3$} \put(4,2){$\oplus$} \put(8,10){$2$}
\put(8,6){$3$}
\put(8,2){$1$}\put(8,-2){$2$}\put(12,2){$\oplus$}\put(16,10){$3$}
\put(16,6){$1$} \put(16,2){$2$}\put(16,-2){$3$}\put(20,-2){$,$}
\end{picture}
\end{center}
\vspace{0.2cm}then there is only two simple-minded systems over
$B'$:

$\mathcal{S}_1' = \{1,2,3\}$; $\mathcal{S}_2' = \{1,\begin{array}{cc} 2 \\
 3
\end{array},\begin{array}{ccc} 3 \\
 1 \\ 2
\end{array}\}$.

This reflects the fact that there are many more (non-trivial) stable
equivalences related to $B$ than that related to $B'$. However, if
we consider the number of orbits of the simple-minded systems under
stable self-equivalences, then in both cases, the number is $2$.
\end{Ex}


\begin{thebibliography}{88}

\bibitem{Asashiba2003}{{\sc H.Asashiba,} On a lift of an individual stable equivalence to a standard
derived equivalence for representation-finite selfinjective
algebras. Algebras and Representation Theory \textbf{6}(4) (2003),
427-447.}

\bibitem{AR1973}{{\sc M.Auslander and I.Reiten,} Stable
equivalence of Artin algebras. In: {\it Proceedings of the
Conference on Orders, Group Rings and Related Topics, Ohio State
Univ., Columbus, Ohio, 1972}. Lecture Notes in Mathematics 353,
Springer, Berlin, 1973, 8-71.}


\bibitem{AR1978}{{\sc M.Auslander and I.Reiten,} Representation theory of Artin algebras
VI, A functorial approach to almost split sequences. Comm. in
Algebra \textbf{6}(3) (1978), 257-300.}



\bibitem{ARS}{{\sc
M.Auslander, I.Reiten and S.O.Smal\o,} {\it Representation theory of
Artin algebras}. Cambridge University Press, 1995.}


\bibitem{Broue1994}{{\sc M.Brou\'e,} Equivalences of blocks of group
    algebras. In: {\it Finite
dimensional algebras and related topics}. V.Dlab and L.L.Scott
(eds.), Kluwer, 1994, 1-26.}

\bibitem{Carlson1998}{{\sc J.Carlson,} A characterization of endotrivial modules over p-groups. Manuscripta
Math. \textbf{97} (1998), 303-307.}

\bibitem{DM2007}{{\sc A.S.Dugas and R.Martinez-Villa,} A note on stable equivalence of Morita
type. J. Pure Appl. Algebra \textbf{208}(2) (2007), 421-433.}

\bibitem{EK2000}{{\sc K.Erdmann and O.Kerner,} On the
stable module category of a self-injective algebra. Trans. Amer.
Math. Soc. \textbf{352} (2006), 2389-2405.}

\bibitem{KonigLiu2008}{{\sc
S.Koenig and Y.M.Liu,} Gluing of idempotents, radical embeddings
and two classes of stable equivalences. J. Algebra \textbf{319}
(2008), 5144-5164.}

\bibitem{Linckelmann1996}{{\sc M.Linckelmann,}
Stable equivalences of Morita type for self-injective algebras and
p-groups. Math. Zeit. \textbf{223} (1996), 87-100.}

\bibitem{Liu2003}{{\sc Y.M.Liu,} On stable
equivalences of Morita type for finite dimensional algebras. Proc.
Amer. Math. Soc. \textbf{131} (2003), 2657-2662.}

\bibitem{Liu2006}{{\sc
Y.M.Liu,} On stable equivalences induced by exact functors. Proc.
Amer. Math. Soc. \textbf{134} (2006), 1605-1613.}

\bibitem{Liu2008}{{\sc Y.M.Liu,} Summands of
stable equivalences of Morita type. Comm. in Algebra \textbf{36}(10)
(2008), 3778-3782.}

\bibitem{LiuXi2006}{{\sc Y.M.Liu and C.C.Xi,}
Constructions of stable equivalences of Morita type for finite
dimensional algebras I. Trans. Amer. Math. Soc. \textbf{358} (2006),
2537-2560.}



\bibitem{MV1981}{{\sc
R.Martinez-Villa,} Algebras stably equivalent to factors of
hereditary. In: {\it Representations of algebras, Puebla, 1980}.
Lecture Notes in Mathematics 903, Springer, Berlin-New York, 1981,
222-241.}



\bibitem{MV1990}{{\sc
R.Martinez-Villa,} Properties that are left invariant under stable
equivalence. Comm. in Algebra \textbf{18}(12) (1990), 4141-4169.}

\bibitem{MV1991}{{\sc
R.Martinez-Villa,} The stable group of a selfinjective Nakayama
algebra. Comm. in Algebra \textbf{19}(2) (1991), 509-517.}

\bibitem{Pogorzaly1993}{{\sc
Z.Pogorza\l y,} On the stable Grothendieck groups. In: {\it Canadian
Mathematical Society Conference Proceedings Volume 14}. 1993,
393-406. }

\bibitem{Pogorzaly1994}{{\sc Z.Pogorza\l
y,} Algebras stably equivalent to self-injective special biserial
algebras. Comm. in Algebra \textbf{22}(4) (1994), 1127-1160.}





\bibitem{Rickard2002}{{\sc J.Rickard,} Equivalences of derived categories for symmetric algebras. J. Algebra  \textbf{257}  (2002), 460-481.}

\bibitem{RR2010}{{\sc J.Rickard and R.Rouquier,} Stable categories and reconstruction. arXiv:1008.1976v1  (2010), 1-16.}

\bibitem{Rouquier2008}{{\sc R.Rouquier,} Dimensions of triangulated categories. Journal of K-theory \textbf{1} (2008), 193-256 and errata, 257-258.}





\end{thebibliography}
\end{document}